\documentclass[12pt]{amsart}
\usepackage[all]{xy} 
\usepackage{amsmath}
\usepackage{amsfonts}
\usepackage{amssymb}

\def\yb{\overline{y}}
\def\jbar{\overline{\jmath}}
\def\etabar{\overline{\eta}}
\def\deltabar{\overline{\delta}}
\def\ra{\rightarrow}
\def\e{\kern 0.08em}
\def\le{\kern 0.05em}
\def\ng{\kern -0.03em}
\def\be{\kern -0.015em}
\def\lbe{\kern -0.01em}
\def\g{\varGamma}
\def\kb{\overline{k}}
\def\cb{\overline{C}}
\def\ub{\overline{U}}
\def\krn{{\rm{Ker}}\,}
\def\kc{\kb\be\!\phantom{.}^{*}}
\def\img{{\rm{Im}}\,}
\def\cok{{\rm{Coker}}\,}

\def\sb{\overline{S}}
\def\xb{\overline{X}}

\def\kc{\kb\!\phantom{.}^{*}}
\def\dimn{\text{dim}\,}
\def\spec{{\rm{Spec}}\,}
\def\bp{{\mathbb P}}

\def\bz{{\mathbb Z}}
\def\bq{{\mathbb Q}}
\def\s{\mathcal }

\newtheorem{lemma}{Lemma}[section]
\newtheorem{theorem}[lemma]{Theorem}
\newtheorem{corollary}[lemma]{Corollary}
\newtheorem{proposition}[lemma]{Proposition}
\theoremstyle{definition}
\newtheorem{definition}[lemma]{Definition}
\theoremstyle{remark}
\newtheorem{remark}[lemma]{Remark}
\newtheorem{remarks}[lemma]{Remarks}

\begin{document}

\title[Algebraic cycles on quadric fibrations]
{Finiteness theorems for algebraic cycles of small codimension on
quadric fibrations over curves}


\author{Cristian D. Gonz\'alez-Avil\'es}
\address{Departamento de Matem\'aticas, Universidad de La
Serena, La Serena, Chile}

\email{cgonzalez@userena.cl}


\thanks{The author is partially supported by Fondecyt grant
1080025}

\maketitle

\begin{abstract} We obtain finiteness theorems for algebraic cycles
of small codimension on quadric fibrations over curves over perfect
fields. For example, if $k$ is finitely generated over $\bq$ and
$X\ra C$ is a quadric fibration of odd relative dimension at least
$11$, then $CH^{\e i}(X)$ is finitely generated for $i\leq 4$.
\end{abstract}

\section{Introduction.}

A well-known conjecture of S.Bloch asserts that the Chow ring of a
smooth projective variety over a number field is a finitely
generated abelian group. In connection with this conjecture, a
number of authors have studied 0-cycles (i.e., cycles of maximal
codimension) on quadric fibrations $\pi\colon X\ra C$ of relative
dimension $d\geq 1$ over smooth integral curves $C$. In \cite{G},
M.Gros studied 0-cycles of degree 0 on conic fibrations (i.e.,
$d=1$) over a number field $k$. The main result of that paper,
obtained by $K$-theoretic methods, was the finiteness of
$\krn(A_{0}(X)\ra A_{0}(\xb)^{\g})$, where $A_{0}(X)$ denotes the
Chow group of 0-cycles of degree 0 on $X$, $\xb=X\times_{\spec
k}\spec \kb$ and $\g=\text{Gal}\big(\e\kb/k\big)$. Further progress
was made by J.-L.Colliot-Th\'el\`ene and A.Skorobogatov in
\cite{CTS}. These authors established the finiteness of the group
$$
CH_{0}(X/C)=\krn[CH_{0}(X)\overset{\pi_{*}}\longrightarrow
CH_{0}(C)]
$$
when $k$ is a number field or a local field and $d=2$. To obtain
this result, they first established an isomorphism
$$
CH_{0}(X/C)=k(C)^{*}_{\text{dn}}/k[C]^{*}N_{X_{\eta}}(k(C))
$$
for any fibration $X\ra C$ satisfying certain assumptions, chiefly
the vanishing of the groups $A_{0}(X_{y})$ associated to the closed
fibers $X_{y}$. Here $k(C)^{*}_{\text{dn}}$ is a certain subgroup of
$k(C)^{*}$ of ``divisorial norms" and $N_{X_{\eta}}\be(k(C))$ is the
group of norms associated to the generic fiber $X_{\eta}$. Then they
used the above isomorphism to reduce the study of the group
$CH_{0}(X/C)$ for a quadric fibration $X\ra C$ of relative dimension
2 to that of $CH_{0}(Y/\widetilde{C}\e)$, where $Y\ra\widetilde{C}$
is a certain {\it conic} fibration over a ``discriminant curve"
$\widetilde{C}$ covering $C$. Thus, using the results of Gros
mentioned above, these authors were able to establish the finiteness
of $CH_{0}(X/C)$ when $\widetilde{C}$ is geometrically integral and
$k$ is a number field or a local field. The next step in the study
of $CH_{0}(X/C)$ was taken by R.Parimala and V.Suresh \cite{PS}, who
developed the methods of \cite{CTS} further. These authors
established the vanishing of $A_{0}(X)$ for {\it smooth Pfister}
quadric fibrations $X$ over {\it conics} $C$ (of arbitrary relative
dimension and defined over any field of characteristic different
from 2). They showed, further, that $A_{0}(X)$ vanishes as well if
$X_{\eta}$ is arbitrary and either $C=\bp^{1}$ or $k$ is a number
field or a field of 2-cohomological dimension $\leq 2$. For possibly
non-smooth quadric fibrations over curves $C$ of arbitrary genus,
Parimala and Suresh established the finiteness of $CH_{0}(X/C)$ over
a local field $k$. Over number fields, they were able to establish
the finiteness of $CH_{0}(X/C)$ (and therefore the finite generation
of $CH_{0}(X)$) for certain ``admissible" quadric fibrations whose
generic fiber is defined by a Pfister neighbor of dimension at least
5. Such is the present state of progress towards obtaining a proof
of Bloch's conjecture for quadric fibrations over curves. In
particular, until now only 0-cycles (i.e., cycles of maximal
codimension) on such fibrations had been studied. In this paper we
study cycles of small codimension on quadric fibrations as above.
Using methods analogous to those developed in \cite{GA}, we obtain
the following result.

\begin{theorem} Let $k$ be a perfect field of characteristic
different from 2 and let $C$ be a smooth, projective and
geometrically integral $k$-curve. Let $X\ra C$ be a quadric
fibration of relative dimension $d\geq 11$, where $X$ is a smooth,
projective and geometrically integral $k$-variety. If $d$ is even,
assume that ${\rm{disc}}\e(X_{\etabar})=1$. Assume, in addition,
that one of the following conditions holds:
\begin{enumerate}
\item[(a)] $k$ is finitely generated over $\,\bq$, or
\item[(b)] $C$ is a conic.
\end{enumerate}
Then $CH^{\e i}(X)$ is finitely generated for $i\leq 4$.
\end{theorem}

The methods of this paper also yield finiteness results for cycles
of codimension $i$ for every $i$ in a certain extended range if the
generic fiber of $X\ra C$ is an excellent quadric. We illustrate
this fact in Section 5 by considering Pfister quadric fibrations.

\section*{Acknowledgements}
I thank R.de Jeu for sending me a proof of Lemma 3.7. I also thank
N.Karpenko for some helpful comments and for sending me copies of
\cite{K1} and \cite{K5} . Finally, I thank B.Kahn, T.Szamuely,
B.Totaro and a number of participants of the seminar Vari\'et\'es
rationnelles (E.N.S., Paris) for additional helpful comments.

\section{Preliminaries.}

Let $k$ be a perfect field of characteristic different from 2 and
let $\kb$ be a fixed algebraic closure of $k$. For any $k$-variety
$Y$, $Y_{0}$ will denote the set of closed points of $Y$. If $Y$ is
a smooth, projective and geometrically integral quadric over $k$,
$\text{disc}(Y)\in k^{*}/(k^{*})^{2}$ will denote the discriminant
(signed determinant) of any quadratic form defining $Y$ \cite{L},
p.38. Now let $\g=\text{Gal}\big(\e\kb/k\big)$ and let $C$ be a
smooth, projective and geometrically integral $k$-curve with
function field $k(C)$ and generic point $\eta$.

\begin{definition} An {\it admissible quadric fibration over $C$}, of relative
dimension $d\geq 1$, is a pair $(X,\pi)$ consisting of a {\it
smooth}, projective and geometrically integral $k$-variety $X$ and a
morphism $\pi\colon X\ra C$ such that each point $y\in C$ has an
affine neighborhood $\spec A(y)$ with $X\times_{C}\,\spec A(y)$
isomorphic to a (possibly singular) quadric in $\bp_{\!A(y)}^{\e
d+1}$. Further, $X_{\eta}=X\times_{C}\spec k(C)$ is smooth and
$\text{disc}(X_{\etabar})=1$ if $d$ is even, where
$X_{\etabar}=X\times_{C}\spec \kb(C)$.
\end{definition}

An admissible quadric fibration $(X,\pi)$ as above will often be
denoted by $X\ra C$.

\begin{remark} Since we have assumed that $X$ is smooth, the class
of admissible quadric fibrations considered in this paper is
narrower than that considered in \cite{PS}. This smoothness
condition is imposed in order to have available the localization
exact sequence \eqref{loc} below. Note, however, that the smooth
Pfister quadric fibrations considered in \cite{PS} are admissible
quadric fibrations in the above sense.
\end{remark}

For any $y\in C$, let $q_{\e y}$ be a fixed quadratic form over
$k(y)$ defining the quadric $X_{y}=X\times_{C}\e\spec k(y)$. We will
write $q_{\e y}^{\e{\rm{ns}}}$ for the nonsingular part of $q_{\e
y}$ and $X_{y}^{{\rm{ns}}}$ for the smooth, projective and
geometrically integral $k(y)$-quadric defined by $q_{\e
y}^{\e{\rm{ns}}}$. Further, set
$$
d_{\le y}=\dimn X_{y}^{{\rm{ns}}}.
$$
Note that, since $X_{\eta}$ is smooth, there exists a finite set $S$
of closed points of $C$ such that, for every $y\in U\!:=C\,\setminus
S$, we have $X_{y}=X_{y}^{{\rm{ns}}}$. Now set $\xb=X\times_{\spec
k}\spec\kb$ and $\cb=C\times_{\spec k}\spec\kb$. The finite set of
closed points of $\cb$ lying above the points in $S$ will be denoted
by $\sb$. Further, $\kb\e[\e U\e]$ will denote the ring of regular
functions on $\ub:=\cb\setminus\sb$. Recall that a smooth,
projective and geometrically integral quadric $Y$ of dimension $d$
over a field $F$ of characteristic not equal to 2 is called {\it
split} if it is isomorphic to either $x_{0}x_{1}+\dots
+x_{d}x_{d+1}=0$ if $d$ is even or $x_{0}^{2}+x_{1}+x_{2}+\dots
x_{d}x_{d+1}=0$ if $d$ is odd. Clearly, if $F$ is algebraically
closed, then any quadric $Y$ over $F$ (as above) is split. On the
other hand, if $F$ is {\it quasi-algebraically closed} (i.e.,
$C_{1}$), then $Y$ is split if, and only if, either $\dimn Y$ is odd
or $\dimn Y$ is even and $\text{disc}(Y)=1$. Thus, if $X\ra C$ is an
admissible quadric fibration, then $X_{\yb}=X\times_{C}\spec
\kb(\yb)$ is a split $\kb(\yb\e)$-quadric for every $\yb\in\ub$
(recall that $\kb(\etabar)=\kb(C)$ is a $C_{1}$-field by Tsen's
theorem).

\begin{lemma} Let $X\ra C$ be an admissible quadric fibration of
relative dimension $d$ and let $i$ be an integer such that $0\leq
i\leq d$.
\begin{enumerate}
\item[(a)] For every $\yb\in\cb$, there exist isomorphisms of abelian groups
$$
CH^{\e i}(X_{\yb})=\begin{cases}
\quad\bz\!\e\qquad\text{if $\,i\neq d_{\e y}/2$}\\
\bz\oplus\bz\quad\text{if $\,i=d_{\e y}/2$.}
\end{cases}
$$
If $i\neq d_{\e y}/2$, then $\g$ acts trivially on $CH^{\e
i}(X_{\yb})$.

\item[(b)] For every $y\in C$,
$$
\krn\!\!\left[\e CH^{\e
i}\!\left(X_{y}\right)\overset{{\rm{res}}_{y}}\longrightarrow CH^{\e
i}\!\left(\e X_{\yb}\right)\right]=CH^{\e i}(X_{y})_{{\rm{tors}}}.
$$
\item[(c)] If $\,i<d/2$, then $\,{\rm{res}}_{\eta}\colon
CH^{\e i}(X_{\eta})\ra CH^{\e i}\be\big(\e X_{\etabar}\e\big)$ is
surjective.
\item[(d)] If $i\neq d/2$, the natural map
$$
\kb(C)^{*}=\kb(C)^{*}\otimes CH^{\e i}(X_{\etabar})\ra H^{\e
i}(X_{\etabar},\s K_{i+1}),
$$
induced by the Brown-Gersten-Quillen spectral sequence and cup
product, is an isomorphism of $\g$-modules.

\end{enumerate}
\end{lemma}
\begin{proof} If $\yb\in\ub$, then $X_{\yb}$ is a split quadric and (a) follows
directly from \cite{K1}, \S2.1. If $\yb\in\sb$, then
\begin{equation}\label{bf}
CH^{\e i}(X_{\yb})=\begin{cases} \, CH^{\e i}(X_{\yb}^{{\rm{ns}}})
\quad\text{if $i\leq d_{y}$}\\
\qquad\bz\qquad\quad\text{if $i>d_{y}$}
\end{cases}
\end{equation}
by \cite{K5}, \S1.1, and (a) again follows from \cite{K1}, \S2.1. As
regards (b), if $X_{y}$ is smooth, i.e., $y\in U$, then (b) follows
from \cite{K1}, (2.7). On the other hand, if $y\in S$, then
\eqref{bf} over $k$ and over $\kb$ together with \cite{K1}, \S2.1,
applied to $X_{y}^{{\rm{ns}}}$ show that
$$
\krn\!\!\left[\e CH^{\e
i}\!\left(X_{y}\right)\overset{{\rm{res}}_{y}}\longrightarrow CH^{\e
i}\!\left(\e X_{\yb}\right)\right]=\begin{cases} \, CH^{\e
i}(X_{y}^{{\rm{ns}}})_{\text{tors}}
\quad\text{if $i\leq d_{y}$}\\
\qquad 0\kern 4.5em\text{if $i>d_{y}$}.
\end{cases}
$$
The latter equals $CH^{\e i}(X_{y})_{{\rm{tors}}}$ by the analogue
of \eqref{bf} over $k$, whence (b) follows. Assertion (c) follows
from \cite{K1}, (2.7). Finally, (d) is a particular case of
\cite{KRS}, Proposition 2.2(b).
\end{proof}

Let $X\ra C$ and $i$ be as in the lemma and assume that $i\neq d/2$.
Then there exists a canonical commutative diagram
$$
\xymatrix{k(C)^{*}\otimes CH^{\e
i}(X_{\eta})\ar[d]_{\text{Id.}\otimes\,\text{res}_{\eta}}\ar[r]&
H^{\e i}\!\left(X_{\eta},\e\s K_{i+1}\right)
\ar[d]^{\text{res}}\\
(\kb(C)^{*}\otimes CH^{\e i}(X_{\etabar}))^{\g}\ar[r]^(.55){\sim}&
H^{\e i}\!\left(X_{\etabar},\e\s K_{i+1}\right)^{\e \g},\\
}
$$
where the horizontal maps are induced by the Brown-Gersten-Quillen
spectral sequence and cup-product. The bottom map is an isomorphism
by part (d) of the lemma. Further, since $CH^{\e
i}(X_{\etabar})=\bz$ with trivial $\g$-action by Lemma 2.3(a), we
have
$$
(\kb(C)^{*}\otimes CH^{\e i}(X_{\etabar}))^{\g}=k(C)^{*}\otimes
CH^{\e i}(X_{\etabar})=k(C)^{*}.
$$
Thus, there exists a map
\begin{equation}\label{ro}
\rho_{\e i}\colon H^{\e i}\!\left(X_{\eta},\e\s K_{i+1}\right)\ra
k(C)^{*},
\end{equation}
a canonical isomorphism
$$
\cok\rho_{\e i}\simeq \cok\!\!\left[H^{\e
i}\!\left(X_{\eta},\e{\mathcal
K}_{i+1}\right)\overset{{\rm{res}}}\longrightarrow H^{\e
i}\!\left(X_{\etabar},\e{\mathcal K}_{i+1}\right)^{\g}\e\right]
$$
and an exact sequence
\begin{equation}\label{es1}
k(C)^{*}\!\otimes\cok\!\!\left[\e CH^{\e
i}(X_{\eta})\overset{\text{res}_{\eta}}\longrightarrow CH^{\e
i}\be\big(\e X_{\etabar}\e\big)\right]\ra\cok\rho_{\e i}\ra 0.
\end{equation}
In particular, Lemma 2.3(c) yields the following.

\begin{proposition} If $\,0\leq i<d/2$, then the map $
\rho_{\e i}\colon H^{\e i}\!\left(X_{\eta},\e\s K_{i+1}\right)\ra
k(C)^{*}$ defined above is {\rm surjective}.\qed
\end{proposition}

\begin{remark} Let $w$ be the Witt index of $q_{\eta}$ and assume that
$w<d/2$. Then, if $d/2<i\leq d-w$, there exists a canonical
isomorphism
$$
\cok\!\!\left[CH^{\e
i}(X_{\eta})\overset{{\rm{res}}_{\eta}}\longrightarrow CH^{\e
i}\be\big(\e X_{\etabar}\e\big)\right]=\bz/2.
$$
See \cite{KRS}, Proposition 1.1(c). Thus \eqref{es1} shows that
$\cok\rho_{\e i}$ is a (possibly nontrivial) quotient of
$k(C)^{*}/\{k(C)^{*}\}^{2}$.
\end{remark}

\medskip

We will need the following basic fact: if $M$ is a $\g$-module which
is finitely generated as an abelian group, then $H^{\e 1}(\g,M)$ is
a group of finite exponent. Indeed, if $L$ is a finite Galois
extension of $k$ contained in $\kb$ which trivializes $M$,
$G=\text{Gal}(L/k)$ and $H=\text{Gal}(\e\kb/L)$, then the
inflation-restriction exact sequence in Galois cohomology yields an
exact sequence
$$
0\ra H^{1}(G,M)\ra H^{1}(\g,M)\ra \text{Hom}(H,M_{\text{tors}}).
$$
The left-hand group above is finite by \cite{W}, Corollary 6.5.10,
p.180, and the right-hand group is annihilated by the order of
$M_{\text{tors}}$. This proves our claim. Note that, if $M$ is {\it
free}, then $H^{1}(\g,M)$ is in fact finite.

\section{The basic exact sequence.}

Let $X\ra C$ be an admissible quadric fibration (see Definition 2.1)
of relative dimension $d$ and let $i$ be an integer such that $0\leq
i\leq d$ and $i\neq d/2$. There exists a canonical homomorphism of
$\g$-modules
\begin{equation}\label{deltabar}
\deltabar_{i}\,\colon H^{\e i}(X_{\etabar}, \s
K_{i+1})\ra\displaystyle\bigoplus_{\yb\e\in\e \cb_{0}}CH^{\e
i}(X_{\yb})
\end{equation}
defined as the limit over all nonempty open subsets $\overline{V}$
of $\cb$ of the boundary maps $H^{\e i}(X_{\overline{V}}, \s
K_{i+1})\ra CH^{\e i}(X_{\e\cb\,\setminus\overline{V}})$ arising
from the localization sequence for the triple
$(X_{\overline{V}},\xb,X_{\cb\,\setminus\overline{V}})$ (for a
description of the latter maps, see \cite{R}, (3.7) and (2.1.0) with
$M=K_{*}^{\prime}$ there). Let $\deltabar_{i\e,\e\sb}$ be the
composite
\begin{equation}\label{es2}
H^{\e i}(X_{\etabar},\s
K_{i+1})\overset{\deltabar_{i}}\longrightarrow\displaystyle
\bigoplus_{\yb\e\in\e \cb_{0}}CH^{\e
i}(X_{\yb})\ra\displaystyle\bigoplus_{\yb\e\notin\e\sb\e\cup\e\{\etabar\}}CH^{\e
i}(X_{\yb}),
\end{equation}
where the second map is the canonical projection. Since $i\neq d/2$,
$H^{\e i}(X_{\etabar},\s K_{i+1})$ is canonically isomorphic to
$\kb(C)^{*}$ by Lemma 2.3(d). Further, if
$\yb\e\notin\e\sb\e\cup\e\{\etabar\}$, then $d_{y}=d$ and Lemma
2.3(a) shows that $CH^{\e i}(X_{\yb})=\bz$. Thus
$\deltabar_{i\e,\e\sb}$ may be identified with a map
$\kb(C)^{*}\ra\bigoplus_{\,\yb\e\notin\e\sb\e\cup\e\{\etabar\}}\bz$,
and the description of the map $\deltabar_{i}$ alluded to above (see
\cite{R}, (3.7) and (2.10)) shows that the latter map coincides with
the canonical divisor map
$f\mapsto(\text{ord}_{\e\yb}(f))_{\e\yb\e\notin\e\sb\e\cup\e\{\etabar\}}$.
Consequently, there exist canonical isomorphisms of $\g$-modules
$$
\krn\deltabar_{i\e,\e \sb}=\kb\e[\e U\e]^{*}
$$
and
$$
\cok\deltabar_{i\e,\e \sb}=CH_{0}(\e\ub\e).
$$
Thus the kernel-cokernel exact sequence \cite{M}, Proposition
I.0.24, p.19, associated to \eqref{es2} yields an exact sequence
\begin{equation}\label{kc}
\begin{array}{rcl}
0 \ra \krn\deltabar_{i}\ra\kb\e[\e U\e]^{*}\ra
\displaystyle\bigoplus_{\yb\e\in\e\sb}CH^{\e
i}(X_{\yb})&\ra&\cok\deltabar_{i}\\
&\ra&CH_{0}(\e\ub\e)\ra 0.
\end{array}
\end{equation}

\begin{proposition} Let $\deltabar_{i}$ be the map \eqref{deltabar},
where $0\leq i\leq d$ and $i\neq d/2$. Then $H^{\e
1}\!\left(\g,\krn\deltabar_{i}\e\right)$ is finite.
\end{proposition}
\begin{proof} By \eqref{kc} and Lemma 2.3(a), there exists an exact sequence
of $\g$-modules
$$
0\ra \krn\deltabar_{ i}\ra\kb\e[\e U\e]^{*}\ra A\ra 0,
$$
where $A$ is free and finitely generated. We conclude that there
exists an exact sequence
$$
A^{\g}\ra H^{\e 1}\!\left(\g,\krn\deltabar_{i}\e\right)\ra H^{\e
1}\be\big(\g,\kb\e[\e U\e]^{*}\big).
$$
The image of the left-hand map above is finite since $A^{\g}$ is
finitely generated and $H^{\e
1}\!\left(\g,\krn\deltabar_{i}\e\right)$ is torsion. On the other
hand, by Hilbert's Theorem 90, $H^{\e 1}\be\big(\g,\kb\e[\e
U\e]^{*})$ injects into $H^{\e 1}\be\big(\g,\kb\e[\e U\e]^{*}\be/\e
\kc)$, which is finite since $\kb\e[\e U\e]^{*}\!\big/\,\kc$ is free
and finitely generated (this is a general fact, but in the present
case it suffices to note that $\kb\e[\e U\e]^{*}\!\big/\,\kc$
injects into $\text{Div}_{\sb}\e\big(\e\cb\e\big)$, the group of
divisors on $\cb$ with support in the finite set $\sb\,$). This
completes the proof.
\end{proof}

\begin{remarks} (a) The proof of the proposition shows that there
exists an injection $(\krn\deltabar_{ i})^{\g}\hookrightarrow k[\e
U\e]^{*}$ whose cokernel is finitely generated.

(b) By \eqref{kc}, if $0\leq i\leq d$, $i\neq d/2$ and $S=\emptyset$
(i.e., $\pi\colon X\ra C$ is smooth), then
$$
\krn\deltabar_{i}=\kb\e[\e C\e]^{*}=\kc.
$$
In this case, therefore, $H^{\e
1}\!\left(\g,\krn\deltabar_{i}\e\right)=H^{\e
1}\!\left(\g,\kc\right)=0$ by Hilbert's Theorem 90.
\end{remarks}

Now let $j\colon X_{\eta}\ra X$ and $\jbar\e\colon
X_{\etabar}\ra\xb$ be the canonical embeddings. There exist
canonical exact sequences
$$
H^{\e i}(X_{\eta},{\mathcal
K}_{i+1})\overset{\delta_{i}}{\longrightarrow}\bigoplus_{y\in
C_{0}}CH^{\e i}(X_{y}) \ra CH^{\e i+1}(X)\overset{j^{*}}\to CH^{\e
i+1}(X_{\eta})\ra 0
$$
and
\begin{equation}\label{locseq}
H^{\e i}(X_{\etabar},{\mathcal
K}_{i+1})\overset{\deltabar_{i}}{\longrightarrow}\bigoplus_{\yb\in
\cb_{0}}CH^{\e i}(X_{\yb}) \ra CH^{\e
i+1}(\e\xb\e)\overset{\overline{\jmath}^{\e *}}\to CH^{\e
i+1}(X_{\etabar})\ra 0
\end{equation}
which yield the following commutative diagrams:

\begin{equation}\label{cd1}
\xymatrix{\krn\delta_{i}\,\ar@{^{(}->}[r]\ar[d]& H^{\e
i}\!\left(X_{\eta},\e\s
K_{i+1}\right)\ar[r]\ar[d]^{\text{res}}&\img\delta_{i}\ar[d]^{\Phi_{i}}
\ar[r]&0&\\
(\e\krn\deltabar_{i})^{\e\g}\,\ar@{^{(}->}[r]& H^{\e
i}\!\left(X_{\etabar},\e\s K_{i+1}\right)^{\g}\ar[r] &
\left(\img\deltabar_{i}\right)^{\g}\ar@{->>}[r]&H^{\e
1}\!\left(\g,\krn\deltabar_{i}\e\right),}
\end{equation}
where $\Phi_{i}$ is induced by the restriction map res and the
bottom row is exact by Hilbert's Theorem 90 via Lemma 2.3(d),
\begin{equation}\label{cd2}
\xymatrix{ 0\ar[r]&\img\delta_{i}
\ar[r]\ar[d]^{\Phi_{i}}&\displaystyle\bigoplus_{y\in C_{0}}CH^{\e
i}\!\left(X_{y}\right) \ar[d]^{\bigoplus\be\text{res}_{y}}
\ar@{->>}[r] & {\rm{Ker}}\,j^{*}\ar[d]\\
0\ar[r]& \left(\img\deltabar_{i}\right)^{\g}\ar[r]&
\displaystyle\bigoplus_{y\in C_{0}}CH^{\e i}\!\left(\e
X_{\yb}\right)^{\g_{\ng y}}{\kern -.5em}\ar[r]
&\left({\rm{Ker}}\,\overline{\jmath}^{\,*}\right)^{\g}
\\
}
\end{equation}
and
\begin{equation}\label{cd3}
\xymatrix{0\ar[r]&{\rm{Ker}}\,j^{*}\ar[r]\ar[d]& CH^{\e
i+1}(X)\ar[r]\ar[d]^{\text{res}}&CH^{\e i+1}(X_{\eta})\ar[d]
\ar[r]&0\\
0\ar[r]&\left({\rm{Ker}}\,\overline{\jmath}^{\,
*}\right)^{\g}\ar[r]& CH^{\e i+1}\!\left(\e\xb\e\right)^{\g}\ar[r] &
CH^{\e i+1}\!\left(\e X_{\etabar}\e\right)^{\g}& &
,\\
}
\end{equation}
where, for each $y\in C_{0}$, we have fixed a point $\yb\in\cb_{0}$
lying above $y$ and written $\g_{\ng y}={\rm{Gal}}\left(\e
\kb/k(y)\right)$. Applying the snake lemma to diagram \eqref{cd1}
and identifying $(\e\krn\deltabar_{i})^{\e\g}$ with a subgroup of
$H^{\e i}\!\left(X_{\etabar},\e\s K_{i+1}\right)^{\g}=k(C)^{*}$, we
obtain the following result.

\begin{proposition} There exists a canonical
exact sequence
$$
0\ra k(C)^{*}/\,(\img\rho_{\e i})\be(\e\krn\deltabar_{i})^{\e\g}\ra
\cok\Phi_{i}\ra H^{\e 1}\!\left(\g,\krn\deltabar_{i}\e\right)\ra
0.\qed
$$
\end{proposition}

Now set
\begin{equation}\label{kres}
CH^{\e i+1}(X)^{\e\prime}=\krn\!\!\left[\e CH^{\e
i+1}\!\left(X\right)\overset{\text{res}}\longrightarrow CH^{\e
i+1}(\e\xb\e)^{\be\g}\e\right]
\end{equation}
and let
\begin{equation}\label{psi}
\Psi_{i}\,\colon\cok\Phi_{i}\ra\displaystyle{\bigoplus_{y\in
C_{0}}}\,CH^{\e i}\!\left(\e X_{\yb}\right)^{\g_{\ng
y}}\!\big/\,\text{res}_{\le y}\e CH^{\e i}\!\left(X_{y}\right)
\end{equation}
be induced by the map
$\left(\img\deltabar_{i}\right)^{\g}\ra\bigoplus_{\,y\in
C_{0}}CH^{\e i}\!\left(\e X_{\yb}\right)^{\g_{\ng y}}$ appearing on
the bottom row of diagram \eqref{cd2}. By Lemma 2.3(b), the kernel
of the middle vertical map in \eqref{cd2} is $\bigoplus_{y\in
C_{0}}\! CH^{\e i}(X_{y})_{\text{tors}}$. Thus, applying the snake
lemma to \eqref{cd2} and using \eqref{cd3} together with Lemma
2.3(b) (for $y=\eta$), we obtain

\begin{proposition} There exists a canonical exact sequence
$$
\begin{array}{rcl}
\krn\Phi_{i}\hookrightarrow\displaystyle\bigoplus_{y\in C_{0}}CH^{\e
i}(X_{y})_{{\rm{tors}}}&\!\!\ra\!\!&\krn\!\left[\,CH^{\e
i+1}(X)^{\e\prime}\!\overset{j^{*}}\ra CH^{\e
i+1}(X_{\eta})_{{\rm{tors}}}\e\right]\\
&\!\!\ra\!\!&\krn\Psi_{i}\ra 0,
\end{array}
$$
where $CH^{\e i+1}(X)^{\e\prime}$ and $\Psi_{i}$ are given by
\eqref{kres} and \eqref{psi}, respectively.\qed
\end{proposition}

We now note that the exact sequence of Proposition 3.3 induces an
exact sequence
$$
0\ra\krn\Psi_{i}^{\e\prime}\ra \krn\Psi_{i}\ra H^{\e
1}\!\left(\g,\krn\deltabar_{i}\e\right),
$$
where
$$
\Psi_{i}^{\e\prime}\colon k(C)^{*}/\,(\img\rho_{\e
i})\be(\e\krn\deltabar_{i})^{\e\g}\ra\displaystyle{\bigoplus_{y\in
C_{0}}}\,CH^{\e i}\!\left(\e X_{\yb}\right)^{\g_{\ng
y}}\!\big/\,\text{res}_{\e y}\e CH^{\e i}\!\left(X_{y}\right)
$$
is the composition of \eqref{psi} and the injection
\begin{equation}\label{inj}
k(C)^{*}/\,(\img\rho_{\e
i})(\e\krn\deltabar_{i})^{\e\g}\hookrightarrow \cok(\Phi_{i})
\end{equation}
coming from Proposition 3.2. Thus, if $i\neq d/2$, then Proposition
3.1 shows that $\krn\Psi_{i}$ is finite if, and only if,
$\krn\Psi_{i}^{\e\prime}$ is finite. Now \eqref{inj} is induced by
the composite
$$\begin{array}{rcl}
k(C)^{*}&\underset{\text{2.3(d)}}\simeq& H^{\e
i}\!\left(X_{\etabar},\e\s K_{i+1}\right)^{\g}\overset{\deltabar_{
i}}\longrightarrow\displaystyle{\bigoplus_{y\in C_{0}}}\,CH^{\e
i}\!\left(\e
X_{\yb}\right)^{\g_{\ng y}}\\
&\overset{\text{proj.}}\longrightarrow
&\displaystyle{\bigoplus_{y\in C_{0}}}\,CH^{\e i}\!\left(\e
X_{\yb}\right)^{\g_{\ng y}}\!\big/\,\text{res}_{\e y}\e CH^{\e
i}\!\left(X_{y}\right)
\end{array}
$$
and we define the {\it $i$-th Salberger group of $X\ra C$},
${\rm{Sal}}_{\e i}(X/C)$, to be the kernel of the preceding
composition, i.e.,
\begin{equation}\label{sal}
{\rm{Sal}}_{\e i}(X/C)=\left\{f\in k(C)^{*}\colon\forall\e y\in
C_{0},\,\deltabar_{i,\e y}(f)\in\text{res}_{y}\e CH^{\e
i}(X_{y})\right\},
\end{equation}
where $\deltabar_{i,\e y}$ is the $y$-component of the composition
$$
k(C)^{*}\overset{\sim}\ra H^{\e i}\!\left(X_{\etabar},\e\s
K_{i+1}\right)^{\g}\overset{\deltabar_{i}}\longrightarrow
\displaystyle\bigoplus_{\,y\in C_{0}}\,CH^{\e i}\!\left(\e
X_{\yb}\right)^{\g_{\ng y}}.
$$
Thus
$$
\krn\Psi_{i}^{\e\prime}={\rm{Sal}}_{\e i}(X/C)/\,(\img\rho_{\e
i})\be(\e\krn\deltabar_{i})^{\e\g}.
$$
The preceding discussion and Proposition 3.1 yield the following
result.
\begin{proposition} Assume that $0\leq i\leq d$ and $\,i\neq d/2$. Then there exists
a canonical exact sequence
$$
0\ra{\rm{Sal}}_{\e i}(X/C)/\,(\img\rho_{\e
i})\be(\e\krn\deltabar_{i})^{\e\g}\ra\krn\Psi_{i}\ra H^{\e
1}\!\left(\g,\krn\deltabar_{i}\e\right),
$$
where $\Psi_{i}$ is the map \eqref{psi} and ${\rm{Sal}}_{\e i}(X/C)$
is the group \eqref{sal}. In particular, $\krn\Psi_{i}$ is finite
if, and only if, ${\rm{Sal}}_{\e i}(X/C)/\,(\img\rho_{\e
i})\be(\e\krn\deltabar_{i})^{\e\g}$ is finite.\qed
\end{proposition}

The basic exact sequence alluded to in the heading of this Section
is the following.

\begin{theorem} Let $X\ra C$ be an admissible quadric fibration of relative
dimension $d$ and let $i$ be an integer such that $0\leq i< d/2$.
Then there exists a canonical exact sequence
$$
\displaystyle{\bigoplus_{y\in C_{0}}}\, CH^{\e
i}(X_{y})_{\le{\rm{tors}}}\ra {\rm{Ker}}\!\left[\,CH^{\e
i+1}(X)^{\e\prime}\!\ra\! CH^{\e
i+1}(X_{\eta})_{{\rm{tors}}}\e\right]\ra H^{\e
1}\!\left(\g,\krn\deltabar_{i}\e\right),
$$
where $CH^{\e i+1}(X)^{\e\prime}$ is the group \eqref{kres}.
\end{theorem}
\begin{proof} By Proposition 2.4, $\rho_{\e i}\colon
H^{\e i}(X_{\eta},\mathcal K_{i+1})\ra k(C)^{*}$ is surjective. Thus
$\img\rho_{\e i}={\rm{Sal}}_{\e i}(X/C)=k(C)^{*}$ and therefore
$$
{\rm{Sal}}_{\e i}(X/C)/\,(\img\rho_{\e
i})\be(\e\krn\deltabar_{i})^{\e\g}=0.
$$
The theorem now follows by combining Propositions 3.4 and 3.5.
\end{proof}

\medskip

We conclude this Section by giving a sufficient condition under
which the group $CH^{\e i+1}(\e\xb\e)^{\be\g}$, appearing in
\eqref{kres}, is finitely generated. Let $J_{C}$ be the Jacobian
variety of $C$.

\begin{lemma} If $J_{C}(k)$ is finitely generated, then so also is
${\rm{Pic}}(\e\ub\e)^{\g}$.
\end{lemma}
\begin{proof} The well-known exact sequence $0\ra
J_{C}(k)\ra {\rm{Pic}}(\e\cb\e)^{\g}\ra\bz$ shows that
${\rm{Pic}}(\e\cb\e)^{\g}$ is finitely generated. Now, by \cite{F},
Proposition 1.8, p.21, there exists a canonical exact sequence of
$\g$-modules
$$
0\ra P_{\e\cb,\e\ub}\ra {\rm{Pic}}(\e\cb\e)\ra
{\rm{Pic}}(\e\ub\e)\ra 0,
$$
where $P_{\e\cb,\e\ub}$ is a finitely generated abelian group. The
above exact sequence induces an exact sequence
$$
{\rm{Pic}}(\e\cb\e)^{\g}\ra {\rm{Pic}}(\e\ub\e)^{\g}\ra
H^{1}(\g,P_{\e\cb,\e\ub}),
$$
where the right-hand group has a finite exponent $m$ (say). It
follows that ${\rm{Pic}}(\e\ub\e)^{\g}$ is a quotient of the inverse
image of ${\rm{Pic}}(\e\cb\e)^{\g}$ under the multiplication-by-$m$
map $m\colon {\rm{Pic}}(\e\cb\e)\ra{\rm{Pic}}(\e\cb\e)$. Since
${\rm{Pic}}(\e\cb\e)_{\e m-\text{tors}}=J_{C}(\e\kb\e)_{\e
m-\text{tors}}$ is finite, the proof is complete.
\end{proof}

\begin{proposition} Let $X\ra C$ be an admissible quadric
fibration of relative dimension $d$ and let $i$ be an integer such
that $0\leq i\leq d$ and $i\neq d/2$. If $J_{C}(k)$ is finitely
generated, then so also is $CH^{\e i+1}(\e\xb\e)^{\be\g}$.
\end{proposition}
\begin{proof} By the exactness of the sequence
$$
0\ra\cok\deltabar_{i}\ra CH^{\e i+1}(\e\xb\e)\ra CH^{\e
i+1}(X_{\etabar})\ra 0
$$
(see \eqref{locseq}) and Lemma 2.3(a) (for $\yb=\etabar$), it
suffices to check that $(\cok\deltabar_{i})^{\g}$ is finitely
generated. The exactness of the sequence
$$
\bigoplus_{\yb\in\sb}CH^{\e i}(X_{\yb})\ra \cok\deltabar_{i}\ra
CH_{0}(\e\ub\e)\ra 0
$$
(see \eqref{kc}) together with Lemma 2.3(a) show that
$(\cok\deltabar_{i})^{\g}$ is finitely generated if
$CH_{0}(\e\ub\e)^{\g}={\rm{Pic}}(\e\ub\e)^{\g}$ is finitely
generated. The result is now immediate from the previous lemma.
\end{proof}

\begin{remark} If $C$ is a conic, then
$J_{C}(k)=0$ is (certainly) finitely generated for any field $k$. On
the other hand, if $k$ is finitely generated over its prime
subfield, then $J_{C}(k)$ is finitely generated by \cite{C},
Corollary 7.2.
\end{remark}

\section{Cycles of codimensions 3 and 4.}

The following statement collects together several results of
N.Karpenko.

\begin{theorem} Let $Y$ be a smooth, projective and geometrically integral quadric
of dimension $d$ over a field of characteristic not equal to 2.
\begin{enumerate}
\item[(a)] $CH^{2}(Y)_{\rm{tors}}$ has order at most 2. If $d\geq
7$, then $CH^{2}(Y)$ is torsion-free.
\item[(b)] $CH^{3}(Y)_{\rm{tors}}$ has order at most 2. If $d\geq
11$, then $CH^{3}(Y)$ is torsion-free.
\item[(c)] $CH^{4}(Y)_{\rm{tors}}$ has order at most 4 if $d\geq 7$.
\end{enumerate}
\end{theorem}
\begin{proof} For (a), see \cite{K1}, Theorem 6.1. For (b) and (c), see
\cite{K3} and \cite{K4}.
\end{proof}

\begin{theorem} Let $k$ be a perfect field of characteristic
different from 2 and let $C$ be a smooth, projective and
geometrically integral $k$-curve. Let $X\ra C$ be an admissible
quadric fibration of relative dimension at least 7, where $X$ is a
smooth, projective and geometrically integral $k$-variety. Then
$\krn\!\!\left[\e CH^{3}\!\left(X\right)\ra
CH^{3}(\e\xb\e)^{\be\g}\e\right]$ is finite.
\end{theorem}
\begin{proof} This follows by combining Theorem 3.6 for $i=2$, Proposition
3.1 and Theorem 4.1(a), (b).
\end{proof}

\begin{corollary} Let the hypotheses be as in Theorem 4.2. Assume,
in addition, that at least one of the following conditions holds:
\begin{enumerate}
\item[(a)] $k$ is finitely generated over $\,\bq$, or
\item[(b)] $C$ is a conic.
\end{enumerate}
Then $CH^{3}(X)$ is finitely generated.
\end{corollary}
\begin{proof} This is immediate from Proposition 3.8, Remark 3.9
and the theorem.
\end{proof}

Similar arguments, using Theorem 4.1(b), (c) in place of Theorem
4.1(a), (b), yield the following result.

\begin{theorem} Let $k$ be a perfect field of characteristic
different from 2 and let $C$ be a smooth, projective and
geometrically integral $k$-curve. Let $X\ra C$ be an admissible
quadric fibration of relative dimension at least 11, where $X$ is a
smooth, projective and geometrically integral $k$-variety. Assume,
in addition, that one of the following conditions holds:
\begin{enumerate}
\item[(a)] $k$ is finitely generated over $\,\bq$, or
\item[(b)] $C$ is a conic.
\end{enumerate}
Then $CH^{4}(X)$ is finitely generated.\qed
\end{theorem}

\section{Pfister quadric fibrations.}

Let $r\geq 3$ and let $q$ be an $r$-fold Pfister form over a field
$F$ of characteristic not equal to 2. Let $Y$ be the quadric defined
by $q$ (so $\text{dim}\,Y=2^{\e r}-2$ and $\text{disc}(Y)=1$). Then
Example 7.3 and Corollary 8.2 of \cite{K6} show that $CH^{\e
i}(Y)_{\text{tors}}=0$ for every $i<2^{\e r-2}$ (as already noted by
N.Karpenko in \cite{K2}). Thus Theorem 3.6 yields the following
result.

\begin{theorem} Let $k$ be a perfect field of characteristic
different from 2 and let $C$ be a smooth, projective and
geometrically integral $k$-curve. Let $X\ra C$ be an $r$-fold
Pfister quadric fibration, i.e., $X_{\eta}$ is defined by an
$r$-fold Pfister form over $k(C)$, where $r\geq 3$. Then, for every
integer $i$ such that $0\leq i\leq 2^{\e r-2}-2$, there exists a
canonical exact sequence
$$
\displaystyle{\bigoplus_{y\in S}}\, CH^{\e
i}(X_{y})_{\le{\rm{tors}}}\ra \krn\!\!\left[\e CH^{\e
i+1}\!\left(X\right)\ra CH^{\e i+1}(\e\xb\e)^{\be\g}\e\right]\ra
H^{\e 1}\!\left(\g,\krn\deltabar_{i}\e\right)\!,
$$
where $S$ denotes the set of points of $C$ where $X\ra C$ has a
singular fiber.\qed
\end{theorem}

\begin{corollary}
Let $k$ be a perfect field of characteristic different from 2 and
let $C$ be a smooth, projective and geometrically integral
$k$-curve. Let $X\ra C$ be a {\rm{smooth}} $r$-fold Pfister quadric
fibration, where $r\geq 3$. Assume that
\begin{enumerate}
\item[(a)] $k$ is finitely generated over $\bq$, or
\item[(b)] $C$ is a conic.
\end{enumerate}
Then $CH^{\e i}(X)$ is finitely generated for $i\leq 2^{\e r-2}-1$.
\end{corollary}
\begin{proof} This is immediate from the above theorem, Proposition
3.1, Proposition 3.8 and Remark 3.9.
\end{proof}

\end{document}